\documentclass{amsart}

\usepackage{enumerate}
\usepackage{amsrefs}

\makeatletter
\renewcommand\MR[1]{%
    \relax\ifhmode\unskip\spacefactor3000 \space\fi
    MR\nolinebreak[3]\hspace{.16667em}%
    \MRno#1 \@nil
    \@ifundefined{BackCite}{}{%
    \ifx\BackCite\back@cite .\fi 
    }%
}
\def\MRno MR#1 #2\@nil{\@stripzeroes#1}
\def\@stripzeroes{%
  \@ifnextchar0{\expandafter\@stripzeroes\@gobble}{}}
\makeatother

\newcommand\p{$p$\nobreakdash}  
\newcommand\h[1]{\nobreakdash-\hspace{0pt}} 
\newcommand\A{\ensuremath{\mathcal{A}}}
\newcommand\C{\ensuremath{\mathcal{C}}}
\newcommand\G{\ensuremath{\mathcal{G}}}
\newcommand\Q{\ensuremath{\mathbb{Q}}}
\newcommand\Khat{\widehat{K}}
\DeclareMathOperator\SL{SL}
\DeclareMathOperator\GL{GL}
\newcommand\xbar{\bar{x}}
\newcommand\iso{\cong}
\newcommand\Z{\ensuremath{\mathbb{Z}}}
\newcommand\into{\hookrightarrow}
\def\<#1>{\langle #1\rangle}

\theoremstyle{plain}
\newtheorem{theorem}{Theorem}[section]
\newtheorem{lemma}[theorem]{Lemma}
\newtheorem{proposition}[theorem]{Proposition}
\newtheorem{corollary}[theorem]{Corollary}

\theoremstyle{remark}
\newtheorem{remark}[theorem]{Remark}

\numberwithin{equation}{section}

\begin{document}

\title{Transitivity properties for group actions on buildings}

\author{Peter Abramenko}
\address{Department of Mathematics\\
University of Virginia\\
Charlotessville, VA 22904}
\email{pa8e@virginia.edu}

\author{Kenneth S. Brown}
\address{Department of Mathematics\\ 
Cornell University\\ 
Ithaca, NY 14853}
\email{kbrown@math.cornell.edu}


\date{January 6, 2006}


\begin{abstract}
We study two transitivity properties for group actions on buildings,
called Weyl transitivity and strong transitivity.  Following hints by
Tits, we give examples involving anisotropic algebraic groups to show
that strong transitivity is strictly stronger than Weyl transitivity.
A surprising feature of the examples is that strong transitivity holds
more often than expected.
\end{abstract}

\maketitle

\section*{Introduction}
\label{sec:introduction}

Suppose a group~$G$ acts by type-preserving automorphisms on a
building~$\Delta$.  If \A\ is a $G$\h-invariant system of apartments
for~$\Delta$, then the action of~$G$ on~$\Delta$ is said to be
\emph{strongly transitive} with respect to~\A\ if it is transitive on
pairs~$(\Sigma,C)$ with $\Sigma\in\A$ and $C$ a chamber in~$\Sigma$.
The theory of strongly-transitive actions is important because of its 
close connection with the theory of groups with a BN-pair \cites{%
  brown89:_build,tits74:_build_bn,ronan89:_lectur,scharlau95:_build}.

There is a weaker notion of transitivity, whose definition makes use
of the ``Weyl-group-valued distance function'' $\delta\colon
\C\times\C \to W$, where $\C=\C(\Delta)$ is the set of chambers
of~$\Delta$ and $W$ is the Weyl group of~$\Delta$.  Namely, we say
that the action of $G$ on~$\Delta$ is \emph{Weyl transitive} if, for
each $w\in W$, the action is transitive on the ordered pairs~$C,C'$ of
chambers such that $\delta(C,C')=w$.  This is equivalent to saying
that $G$ is transitive on~$\C$ and that the stabilizer of a given
chamber~$C$ is transitive on the $w$\h-sphere $\{D\in\C :
\delta(C,D)=w\}$ for every $w\in W$.  As with strong transitivity,
there is a group-theoretic formulation of Weyl transitivity.  This
theory is sketched by Tits in~\cite{tits92:_twin_kac_moody}, and a
full account will appear
in~\cite{abramenko06:_approac_build}*{Chapter~6}.  The structure is
something like a BN-pair, but one only has the $B$ (sometimes called a
\emph{Tits subgroup} of~$G$), and not necessarily the~$N$.

If the building~$\Delta$ is spherical, then the theory simplifies
considerably.  First, there is a unique system of apartments, so one
can talk about strong transitivity without specifying~\A.  Secondly,
strong transitivity turns out to be equivalent to Weyl transitivity.
For non-spherical buildings, on the other hand, strong transitivity
implies Weyl transitivity but not conversely.  To the best of our
knowledge, however, there are no explicit examples in the literature
to show that the converse is false.  All we have found is a general
suggestion by Tits \cite{tits92:_twin_kac_moody}*{Section~3.1,
  Example~(b)}, where he describes a source of possible examples of
Weyl-transitive actions that are not strongly transitive with respect
to any apartment system.  He does not phrase this in terms of
transitivity properties, but rather in group-theoretic terms.  In the
terminology we introduced above, Tits describes a way to exhibit pairs
$(G,B)$ such that $B$ is a Tits subgroup of~$G$ that does not come
from a BN-pair.

For his proposed examples, $G$ is the group~$\G(K)$ of rational points
of a simple simply-connected algebraic group~\G\ over a global
field~$K$, and $\Delta$ is the Bruhat--Tits building associated
with~$\G$ and a non-Archimedean completion~$\Khat$ of~$K$.  (The
$K$\h-rank of~\G\ should be strictly less than its $\Khat$\h-rank;
otherwise Bruhat--Tits theory~\cite{bruhat72:_group} would imply that
the action of~$G$ on~$\Delta$ is strongly transitive with respect to a
suitable apartment system.)  But Tits did not actually give any
specific examples of \G, $K$, and~$\Khat$ for which the action is Weyl
transitive but is not strongly transitive with respect to \emph{any}
apartment system.

The main purpose of this note is to carry out Tits's suggestion in
detail in the simplest possible case, where $K$ is the field~\Q\ of
rational numbers and \G\ is the norm~1 group of a quaternion division
algebra~$D$ over~\Q.  The completion $\Khat$ is then the field~$\Q_p$
of \p-adic numbers for some prime~$p$, and we denote by~$\Delta_p$ the
corresponding building.  Note that $D$ splits over~$\Q_p$ for almost
all primes~$p$, and $\Delta_p$ is then the tree associated
to~$\SL_2(\Q_p)$~\cite{serre80:_trees}.  The group $G=\G(\Q)$ of
rational points is the multiplicative group of units in~$D$ of norm~1,
and our result is the following dichotomy:

\begin{theorem}
\label{thr:1}
Let $D$ be a quaternion division algebra over~\Q, and let $G$ be its
norm~1 group.  Then, with the notation above, one of the following
conditions holds.
\begin{enumerate}[\rm (a)]
\item\label{item:1}
$-1$ has a square root in~$D$, and, for almost all primes~$p$, the
  action of $G$ on~$\Delta_p$ is strongly transitive with respect to
  some apartment system.
\item
$-1$ does not have a square root in~$D$, and, for almost all
  primes~$p$, the action of $G$ on~$\Delta_p$ is Weyl transitive but
  is not strongly transitive with respect to any apartment system.
\end{enumerate}
\end{theorem}

This is a consequence of a more precise result, stated as
Theorem~\ref{thr:2} in Section~\ref{sec:examples}.  We were quite
surprised by this result.  More precisely, we were surprised that
strong transitivity occurs as often as it does, given that our
group~\G\ is \Q\h-anisotropic.  We do not yet understand what happens
for more general~\G, but we have some evidence that strong
transitivity is relatively rare in this context, as Tits suggested.
On the other hand, norm~1 groups of quaternion algebras are not the
only anisotropic examples where strong transitivity can occur.

In order to keep this paper as elementary as possible, we will make no
further reference to the theory of algebraic groups.  Instead, we will
simply work directly with quaternion algebras.  Moreover, we make very
little use of the theory of buildings, beyond standard terminology.
Indeed, the buildings in our examples are the trees associated
with~$\SL_2(\Q_p)$, and everything we need about these can be found in
Serre~\cite{serre80:_trees}.

\section{Notation and preliminaries}
\label{sec:notation}

We assume familiarity with the theory of quaternion algebras, for
which we refer to \citelist{\cite{lam05:_introd}
  \cite{o'meara00:_introd} \cite{scharlau85:_quadr_hermit}}.

Fix nonzero rational numbers $\alpha,\beta$, and let $D$ be the
corresponding quaternion algebra over~\Q, which we denote by
$(\alpha,\beta)_\Q\,$.  It is a 4\h-dimensional associative algebra with
basis $e_1,e_2,e_3,e_4$, where $e_1$ is the identity element,
$e_2^2=\alpha$, $e_3^2=\beta$, and $e_2e_3=-e_3e_2=e_4$.  Here
$\alpha$ and~$\beta$ are identified with $\alpha e_1$
and~$\beta e_1$.  Recall that $D$ is a division algebra if and
only if its norm form~$N$ is anisotropic.  Here the norm of
$x=x_1+x_2e_2+x_3e_3+x_4e_4$ is $N(x):= x\xbar = x_1^2 - \alpha x_2^2
- \beta x_3^2 + \alpha\beta x_4^2 \in \Q$, where $\xbar:=
x_1-x_2e_2-x_3e_3-x_4e_4$.

We assume from now on that $\alpha$ and~$\beta$ have been chosen so
that $D$ is a division algebra.  We can assure this, for example, by
taking $\alpha,\beta<0$.  For any prime~$p$, let $D_p$ be the
quaternion algebra $(\alpha,\beta)_{\Q_p}$ over~$\Q_p$ obtained
from~$D$ by extension of scalars.  Let $G$ (resp.~$G_p$) be the
subgroup of~$D^*$ (resp.~$D_p^*$) consisting of elements of norm~1.
In what follows, we will only be interested in primes~$p$ such
that~$D_p$ splits.  Thus $D_p$ is isomorphic to the
algebra~$M_2(\Q_p)$ of $2\times 2$ matrices, and $G_p$ is isomorphic
to~$\SL_2(\Q_p)$.

If $-1$ has a square root in~$D$, then $D$ is isomorphic to a
quaternion algebra $(\gamma,-1)_\Q$ for some $\gamma\in\Q^*$; see
Bourbaki~\cite{bourbaki58:_elemen}*{Section~11.2, proof of
  Proposition~1} or Lam~\cite{lam05:_introd}*{proof of
  Theorem~III.5.1}.  We may therefore assume without loss of
generality that $\beta=-1$ in this case.

Consider a prime $p\neq 2$ such that $\alpha$ and~$\beta$ are \p-adic
units, i.e., $v_p(\alpha)=v_p(\beta)=0$, where $v_p$ is the \p-adic
valuation.  (Note that almost all primes satisfy these conditions.)
It is then well-known that $D_p$ splits.  It will be convenient for us
to have a specific isomorphism $D_p \to M_2(\Q_p)$, for which we will
use the following lemma:

\begin{lemma}
\label{lem:1}
$D$ is isomorphic to a quaternion algebra $(\alpha',\beta)_\Q$ for
  some $\alpha'\in\Q^*$ such that $\alpha'$ is a \p-adic unit and has
  a square root in~$\Q_p$\,.
\end{lemma}

\begin{proof}
We try to replace the basis vector $e_2\in D$ by a suitable linear
combination $e_2':= \lambda e_2 + \mu e_4$ with $\lambda,\mu\in\Q$.
Note that any such~$e_2'$ anti-commutes with~$e_3$ and that $(e_2')^2=
\lambda^2\alpha - \mu^2\alpha\beta =: \alpha'\in\Q$.  We will show
that $\lambda,\mu$ can be chosen so that $\alpha'\in U^2$, where $U$
is the group~$\Z_p^*$ of \p-adic units.  Setting $e_4'=e_2'e_3$, we
will then have a ``quaternion basis'' $1,e_2',e_3,e'_4$, showing that
$D\iso (\alpha',\beta)_\Q$\,.

The expression defining~$\alpha'$ above is a binary quadratic form
in the variables~$\lambda,\mu$.  Since the coefficients are \p-adic
units, this form represents all elements of~$U$.  (This follows, for
instance, from \cite{serre73}*{Corollary~2 in Section~II.2.2 and
  Proposition~4 in Section~IV.1.7}.)  In particular, it represents~1,
so we can find $\lambda,\mu\in\Q_p$ with $\lambda^2\alpha -
\mu^2\alpha\beta = 1$.  Since $U^2$ is an open subset of~$\Q_p$, we
can replace $\lambda,\mu$ by approximations in~\Q\ and still have
$\lambda^2\alpha - \mu^2\alpha\beta \in U^2$.
\end{proof}

In view of the lemma we can and will assume, without loss of
generality, that $\alpha$ has a square root in~$\Q_p$.  (Note that
$\beta$ does not change in Lemma~\ref{lem:1}, so it is still~$-1$ by
our earlier convention if $-1$ has a square root in~$D$.)  We can now
exhibit a specific isomorphism $D_p \to M_2(\Q_p)$, given by
\begin{equation}
\label{eq:1}
  x_1+x_2e_2+x_3e_3+x_4e_4 \mapsto
  \begin{pmatrix}
    x_1+x_2\sqrt{\alpha} & \beta\left(x_3+x_4\sqrt{\alpha}\right)\\[3\jot]
    x_3-x_4\sqrt{\alpha} & x_1-x_2\sqrt{\alpha}
  \end{pmatrix}\;.
\end{equation}

\section{Density lemmas}
\label{sec:density-lemmas}

It is obvious that $D$ is dense in~$D_p$\,, where the latter is
topologized as a 4\h-dimensional $\Q_p$\h-vector space.  It is also
true, but not obvious, that the density persists when one passes to
elements of norm~1:

\begin{lemma}
\label{lem:2}
$G$ is dense in~$G_p$\,.
\end{lemma}

\begin{proof}
This is a special case of the weak approximation theorem
\cite{platonov94:_algeb}*{Chapter~7}, but we will give a direct proof.
The main point is to construct lots of elements of~$G$, which we do by
the following ``normalization'':  Given $x\in D^*$, let
$x'=x\xbar^{-1} = x^2/N(x)$; then $x'$ has norm~1.  Using this
construction, we see that the closure of $G$ in~$G_p$ contains all
elements of the form~$y^2/N(y)$ with $y \in D_p^*\,$.  In particular,
it contains all squares of elements of~$G_p\,$, so the proof will be
complete if we show that $G_p$ is generated by squares.  This follows,
for instance, from the fact that $G_p \iso \SL_2(\Q_p)$; the latter is
generated by strictly triangular matrices, all of which are squares.
\end{proof}

Now let $T$ be the ``torus'' in $G$ consisting of quaternions of the
form $x=x_1 + x_2 e_2$ with $N(x)=1$, and let $T_p$ be the similarly
defined subgroup of~$G_p$.  Under the identification of $G_p$
with~$\SL_2(\Q_p)$ that one gets from~\eqref{eq:1}, $T_p$ is simply
the standard torus, consisting of the diagonal matrices of
determinant~1.

\begin{lemma}
\label{lem:3}
$T$ is dense in $T_p$
\end{lemma}

\begin{proof}
We use the same normalization trick as in the proof of
Lemma~\ref{lem:2}.  Namely, we construct elements of~$T$ by starting
with an arbitrary $x=x_1+x_2 e_2 \in D^*$ and forming $x':=x^2/N(x)$.
Computing the images of such elements~$x'$ in~$\SL_2(\Q_p)$ under the
map in~\eqref{eq:1}, we find that they are the diagonal matrices with
diagonal entries~$\lambda,\lambda^{-1}$, where
\begin{equation}
\label{eq:2}
  \lambda=\frac{x_1+x_2\sqrt{\alpha}}{x_1-x_2\sqrt{\alpha}}
\end{equation}
for some $x_1,x_2\in \Q$ that are not both zero.  [Note that the
  denominator is not zero since, in view of our assumption that $D$ is
  a division algebra, $\alpha$ does not have a square root in~\Q.]
The closure of $T$ in~$G_p$ therefore contains all diagonal matrices
of the same form, where now $x_1,x_2\in\Q_p$ and the numerator and
denominator are assumed to be nonzero.  To complete the proof, we will
show that every $\lambda\in\Q_p^*$ can be expressed in this way.  Given
$\lambda\in\Q_p^*$, let's first try to achieve this with $x_2=1$, i.e.,
we try to solve
\begin{equation}
\label{eq:3}
  \lambda=\frac{x+\sqrt{\alpha}}{x-\sqrt{\alpha}}
\end{equation}
for $x\in\Q_p$ with $x\neq\pm\sqrt{\alpha}\,$.  Formally
solving~\eqref{eq:3} for~$x$, we find
\[
x=\sqrt{\alpha}\,\frac{\lambda+1}{\lambda-1}\;,
\]
so we are done if $\lambda\neq 1$.  But we can take care of
$\lambda=1$ by putting $x_2 = 0$ in~\eqref{eq:2}.
\end{proof}

Finally, we record for ease of reference a simple observation that we
will use when we apply the density lemmas.

\begin{lemma}
\label{lem:6}
Let $H$ be a topological group and $H'$ a dense subgroup.
\begin{enumerate}[\rm (1)]
\item \label{item:7}
If $U$ is an open subgroup of~$H$, then $H'$ maps onto~$H/U$ under the
quotient map $H\to H/U$.
\item \label{item:8}
If $H$ acts transitively on a set $X$ and the stabilizer of some point
is an open subgroup, then the action of $H'$ on~$X$ is transitive.
\item \label{item:9}
If $H$ acts on an arbitrary set~$X$ and the stabilizers are open
subgroups, then the $H'$\h-orbits in~$X$ are the same as the $H$\h-orbits.
\end{enumerate}
\end{lemma}

\begin{proof}
For (\ref{item:7}), observe that every coset $hU$ is a nonempty open
set, so it meets~$H'$.  (\ref{item:8}) is a restatement
of~(\ref{item:7}), and (\ref{item:9}) follows from~(\ref{item:8}).
\end{proof}

\section{Strong transitivity vs.\ Weyl transitivity}
\label{sec:strong-trans-vs}

We now digress to clarify conceptually the difference between strong
transitivity and Weyl transitivity.  We will then be able to give our
main results in the next section.  We assume familiarity with standard
terminology regarding buildings and apartment systems
\cites{brown89:_build,ronan89:_lectur,tits74:_build_bn,scharlau95:_build}.

\begin{lemma}
\label{lem:4}
Strong transitivity (with respect to some apartment system) implies
Weyl transitivity.
\end{lemma}

\begin{proof}
Assume the action is strongly transitive, and choose a fixed pair
$(\Sigma,C)$ as in the definition of strong transitivity.  We will
show that the stabilizer of~$C$ is transitive on the $w$\h-sphere for
each $w\in W$; this implies Weyl transitivity since the action is
already known to be transitive on the chambers.  Given~$w$, there is
a unique chamber $C_w\in\C(\Sigma)$ with $\delta(C,C_w)=w$.  (If we
identify $\Sigma$ with~$\Sigma(W,S)$ in such a way that $C$
corresponds to the fundamental chamber, then $C_w$ is simply~$wC$.)
Let $D$ be an arbitrary chamber of~$\Delta$ with $\delta(C,D)=w$, and
let $\Sigma'$ be an apartment containing $C$ and~$D$.  By strong
transitivity there is an element $g\in G$ that stabilizes~$C$ and maps
$\Sigma'$ to~$\Sigma$.  Then $\delta(C,gC')=\delta(C,C')=w$, so
$gC'=C_w$.  Thus the stabilizer of~$C$ is transitive on the
$w$\h-sphere, as required.
\end{proof}

If one wants to try, conversely, to show that a given Weyl-transitive
action is strongly transitive with respect to some apartment system,
one needs to first construct a suitable apartment system.  This is
easy:

\begin{lemma}
\label{lem:5}
Suppose the action of $G$ on~$\Delta$ is Weyl transitive, and let
$\Sigma$ be an arbitrary apartment (in the complete system of
apartments).  Then the set $G\Sigma := \{g\Sigma : g\in G\}$ is a
system of apartments.
\end{lemma}

\begin{proof}
It suffices to show that any two chambers~$C,C'$ are contained in
some~$g\Sigma$.  By transitivity of $G$ on~$\C(\Delta)$, we may assume
that $C\in\C(\Sigma)$, in which case we can find $C''\in\C(\Sigma)$
with $\delta(C,C'')=\delta(C,C')$.  Weyl transitivity now gives us a
$g\in G$ that stabilizes~$C$ and takes $C''$ to~$C'$.  Hence $C,C'\in
g\Sigma$, as required.
\end{proof}

Combining the two lemmas, we can clarify the relationship between the
two notions of transitivity:

\begin{proposition}
\label{prop:1}
The following conditions are equivalent for a type-preserving action
of a group~$G$ on a building~$\Delta$.
\begin{enumerate}[\rm (i)]
\item \label{item:2} The $G$\h-action on~$\Delta$ is strongly
transitive with respect to some apartment system.
\item \label{item:3}
The $G$\h-action on~$\Delta$ is Weyl transitive, and there is an
apartment~$\Sigma$ (in the complete system of apartments) such that
the stabilizer of~$\Sigma$ acts transitively on~$\C(\Sigma)$.
\end{enumerate}
\end{proposition}

\begin{proof}
The implication (\ref{item:2}) $\implies$ (\ref{item:3}) is immediate
from Lemma~\ref{lem:4} and the definition of strong transitivity.
Conversely, if (\ref{item:3}) holds, then the action is strongly
transitive with respect to $\A = G\Sigma$, which is an apartment
system by Lemma~\ref{lem:5}.
\end{proof}

Finally, we record a simple method for constructing Weyl-transitive
actions.

\begin{proposition}
\label{prop:2}
Let $G$ act Weyl transitively on a building~$\Delta$, and suppose that
$G$ is a topological group and that the stabilizer~$B$ of some chamber
is an open subgroup.  If $G'$ is a dense subgroup of~$G$, then the
action of $G'$ on~$\Delta$ is also Weyl transitive.
\end{proposition}

\begin{proof}
Consider the diagonal action of $G$ on~$\C\times\C$, where
$\C=\C(\Delta)$.  Note that every stabilizer is an open subgroup
of~$G$, being an intersection of two conjugates of~$B$.  To show that
the action of~$G'$ is Weyl transitive, we must show that the
$G'$\h-orbits in~$\C\times\C$ are the sets of the form $\{(C,C') :
\delta(C,C') = w\}$, one for each $w\in W$.  But these are precisely
the $G$\h-orbits by assumption, so the result follows from
Lemma~\ref{lem:6}(\ref{item:9}).
\end{proof}

Notice that the action of~$G$ might well be strongly transitive, but
there is no reason to think that the same is true of the action
of~$G'$.  We are now in a position to give specific examples of this.

\section{The examples}
\label{sec:examples}

We return to the hypotheses and notation of
Section~\ref{sec:notation}.  In particular, $D=(\alpha,\beta)_\Q$ is a
quaternion division algebra, $p$ is an odd prime such that
$v_p(\alpha) = v_p(\beta) = 0$, $G$ is the norm~1 group of~$D$, and
$G_p$ is the norm~1 group of~$D_p$.  Since $G_p \iso \SL_2(\Q_p)$, we
have a BN-pair in~$G_p$ in a well-known way and a tree~$\Delta_p$ on
which $G_p$ acts~\cite{serre80:_trees}.  The action is strongly
transitive with respect to the complete apartment system.  This is
proved in greater generality in \cite{brown89:_build}*{Section~VI.9F},
and the result in the present context can also be found in
Serre~\cite{serre80:_trees}*{p.~72}.

\begin{proposition}
\label{prop:3}
The action of $G$ on~$\Delta_p$ is Weyl transitive.
\end{proposition}

\begin{proof}
Since $G$ is dense in~$G_p$ by Lemma~\ref{lem:2}, and since the $B$ of
the BN-pair in~$G_p$ is an open subgroup, this follows from
Proposition~\ref{prop:2}.
\end{proof}

\begin{remark}
\label{rem:1}
Note that, just from the fact that $G$ is transitive on the chambers,
we get a decomposition of~$G$ as an amalgamated free
product~\cite{serre80:_trees}, as in the better-known case of~$\SL_2$.
The same is true if $G$ is replaced by any of its subgroups that are
dense with respect to the \p-adic topology.
\end{remark}

Recall that, in the action of $\SL_2(\Q_p)$ on~$\Delta_p$\,, there is
an apartment~$\Sigma_0$, which we call the \emph{standard apartment},
whose stabilizer is the monomial group; the diagonal matrices act
on~$\Sigma_0$ as translations, and the non-diagonal monomial matrices
act as reflections.  The translation action of the diagonal
group~$T_p$ is given by a surjective homomorphism $T_p\to\Z$ whose
kernel is~$T_p\cap B$, which is an open subgroup of~$T_p$\,.  In view of
Lemma~\ref{lem:3} and Lemma~\ref{lem:6}(\ref{item:7}), it follows that
all of the translations can be achieved by elements of~$T$.  But, as
we are about to see, one can not in general realize the reflections by
elements of~$G$.

\begin{theorem}
\label{thr:2}
The following conditions are equivalent:
\begin{enumerate}[\rm (i)]
\item \label{item:4}
$-1$ has a square root in~$D$.
\item \label{item:5} 
$G$ contains an element that stabilizes the
  standard apartment~$\Sigma_0$ and acts as a reflection on it.
\item \label{item:6}
The action of $G$ on $\Delta_p$ is strongly transitive with respect to some
apartment system.
\end{enumerate}
\end{theorem}

\begin{proof}
If (\ref{item:4}) holds, then $\beta=-1$ by our convention in
Section~\ref{sec:notation}.  The quaternion~$e_3$ is therefore in $G$
and maps to $\left(\begin{smallmatrix}
  0&-1\\1&\phantom{-}0 \end{smallmatrix}\right) \in\SL_2(\Q_p)$.  This
proves~(\ref{item:5}).  The latter implies~(\ref{item:6}) by
Proposition~\ref{prop:1}, since the dihedral group of type-preserving
automorphisms of an apartment is generated by the translations and any
one reflection.  Finally, suppose (\ref{item:6}) holds.  Then $G$
contains an element~$g$ that stabilizes an apartment~$\Sigma$ and acts
as a reflection on it.  To prove~(\ref{item:4}), it suffices to note
that any such~$g$ satisfies $g^2=-1$.  In case $\Sigma=\Sigma_0$, this
is immediate, since $g$ must map to a matrix of the form
$\left(\begin{smallmatrix}
  0&-\lambda^{-1}\\\lambda&0 \end{smallmatrix}\right) \in\SL_2(\Q_p)$.
The general case now follows from the fact that $\SL_2(\Q_p)$ acts
transitively on the complete apartment system.
\end{proof}

Note that the dichotomy stated as Theorem~\ref{thr:1} in the
introduction is an immediate consequence of Theorem~\ref{thr:2} and
Proposition~\ref{prop:3}.

To get specific examples of actions that are Weyl transitive but not
strongly transitive, we need to choose~$\alpha,\beta$ so that
$-1\notin D^2$.  Now direct calculation shows that $-1\in D^2$ if and
only if the ternary quadratic form $\<\alpha,\beta,-\alpha\beta>$
represents~$-1$.  [Here we use the angle-bracket notation
  $\<\alpha_1,\alpha_2,\dots,\alpha_n>$ for the quadratic form
  $\sum_{i=1}^n \alpha_ix_i^2$ in $n$ variables $x_1,\dots,x_n$.]  Hence
\begin{equation}
\label{eq:4}
  -1 \notin D^2 \iff \<1,\alpha,\beta,-\alpha\beta> \text{ is
    anisotropic.}
\end{equation}
Let $l$ be a prime such that $l\equiv 1$ mod~4, so that $-1\in\Q_l^2$.
Set $\beta=-l$, and let $\alpha$ be any negative integer such that
$\alpha$ is not a square mod~$l$.  Then the quaternary form
in~\eqref{eq:4} is equivalent over~$\Q_l$ to the form
$\<1,\alpha,l,\alpha l>$, which is easily seen to be anisotropic
over~$\Q_l$.  In fact, it is the essentially unique anisotropic
quaternary form over~$\Q_l$
\citelist{\cite{lam05:_introd}*{Theorem~VI.2.2(3)}
  \cite{serre73}*{Section~IV.2.3, Corollary to Theorem~7}
  \cite{o'meara00:_introd}*{63:17}}.  The form is therefore
anisotropic over~\Q, so $-1$ does not have a square root in
$D:=(\alpha,-l)_\Q$.  For a concrete example, take $l=5$ and
$\alpha = -2$.

\begin{corollary}
\label{cor:1}
Let $D$ be the quaternion division algebra $(-2,-5)_\Q$\,, and let $G$
be its norm~1 group.  Then for all primes $p\neq 2,5$, there is a
Weyl-transitive action of~$G$ on~$\Delta_p$ that is not strongly transitive
with respect to any apartment system, where $\Delta_p$ is the regular tree
of degree $p+1$.  \qed
\end{corollary}

\section{The role of the apartment system}
\label{sec:role-apartm-syst}

As we have emphasized from the beginning, one needs a system of
apartments in order to talk about strong transitivity.  In our
examples, however, either strong transitivity fails regardless of the
apartment system or strong transitivity holds with respect to the
``standard'' apartment system~$\A_0:=G\Sigma_0$.  This raises the
question of how much choice there is in finding an apartment
system~\A\ such that a given action is strongly transitive with
respect to~\A.  In order to shed some light on this, we consider an
even simpler situation than in the previous section, namely, we take
$G=\SL_2(\Q)$ and consider its natural action on~$\Delta_p$ via the
inclusion $G\into G_p:=\SL_2(\Q_p)$.  The action is strongly
transitive with respect to $\A_0:=G\Sigma_0$, but we will see that it
is also strongly transitive with respect to other apartment
systems.

Fix a prime $p$, and suppose that we have matrices $A,B\in
G=\SL_2(\Q)$ with the following properties:
\begin{itemize}
\item
The characteristic polynomial of~$A$ is irreducible over~\Q\ but
splits over~$\Q_p$.
\item
The eigenvalues $\lambda,\lambda^{-1}$ of~$A$ have \p-adic
valuation~$\pm1$.
\item
$BAB^{-1} = A^{-1}$.
\end{itemize}
Then $A$ is diagonalizable over~$\Q_p$, and, in fact, we can find
$g\in\SL_2(\Q_p)$ such that
\[
gAg^{-1} = \begin{pmatrix}
\lambda&0\\ 0&\lambda^{-1}
\end{pmatrix}\;.
\]
(To see that $g$ can be taken to have determinant~1, observe that the
centralizer in~$\GL_2(\Q_p)$ of the diagonal matrix above is the full
diagonal group, which contains matrices of arbitrary nonzero
determinant.)

It follows easily that the stabilizer of~$\Sigma_0$ in~$gGg^{-1}$ acts
transitively on~$\C(\Sigma_0)$.  Indeed, $gAg^{-1}$
stabilizes~$\Sigma_0$ and generates the infinite cyclic group of
type-preserving translations of the latter; and $gBg^{-1}$ is
necessarily a non-diagonal monomial matrix, since it conjugates
$\left(\begin{smallmatrix}\lambda&0\\ 0&\lambda^{-1}
\end{smallmatrix}\right)$ to $\left(\begin{smallmatrix}\lambda^{-1}&0\\
  0&\lambda\end{smallmatrix}\right)$, so it acts as a reflection
  on~$\Sigma_0$.  

Consequently, the stabilizer of~$\Sigma:=g^{-1}\Sigma_0$ in~$G$ acts
transitively on~$\C(\Sigma)$.  Thus $G$ is strongly transitive with
respect to the apartment system $\A:=G\Sigma$ by
Proposition~\ref{prop:1} (and its proof).  Notice that $\A\neq\A_0$,
i.e., $\Sigma\notin\A_0$, because the matrix~$A$ acts as a translation
on~$\Sigma$; but every element of~$G$ that acts as a translation on an
apartment in~$\A_0$ is diagonalizable over~\Q.

It is easy to find specific examples of the situation we have just
described.  With $p=3$, for instance, we can take
\[
A=\begin{pmatrix}
0&-1\\ 1&7/3
\end{pmatrix} \quad \text{and} \quad
B=\begin{pmatrix}
2&3\\ -5/3&-2
\end{pmatrix}\;.
\]
One of the eigenvalues of~$A$ is $\lambda=(7+\sqrt{13})/6 \in\Q_3$;
here $\sqrt{13}$ exists in~$\Q_3$ because $13\equiv 1$ mod~3, and we
choose the square root that is also $\equiv 1$ mod~3.  Then
$7+\sqrt{13}$ is a 3-adic unit, so $v_3(\lambda) = -1$.  One can check
by direct calculation that $BAB^{-1} = A^{-1}$.

\begin{remark}
\label{rem:2}
To get strong transitivity, one cannot simply use~$G\Sigma$ for an
arbitrary apartment~$\Sigma$ in the complete apartment system.
Indeed, an apartment is completely determined by a single element that
acts as a non-trivial translation on it.  (The element is a
``hyperbolic'' automorphism of the tree, and the apartment is its
axis, cf.~Serre~\cite{serre80:_trees}*{Section~I.6.4}.)  So there are
only countably many choices of~$\Sigma$ that will work.  But the
complete apartment system is uncountable.
\end{remark}

\end{document}